\documentclass[12pt]{amsart}
\usepackage{psfig,graphicx,amsfonts,amssymb,amsmath,amsthm}

\newcommand{\Z}{\mathbb Z} % integers
\newcommand{\Q}{\mathbb Q} % rationals
\newcommand{\R}{\mathbb R} % real
\newcommand{\ZHS}{$\mathbb Z$\rm HS}
\newcommand{\QHS}{$\mathbb Q$\rm HS}

\newcommand{\cB}{\mathcal{B}}
\newcommand{\M}{\mathcal{M}}
\newcommand{\cP}{\mathcal{P}}
\newcommand{\cM}{\mathcal{M}}
\newcommand{\cF}{\mathcal{F}}
\newcommand{\grad}{\mathrm{Gr}}
\newcommand{\cA}{\mathcal{A}}
\newcommand{\frakg}{\mathfrak{g}}
\newcommand{\fg}{\mathfrak{g}}

\newcommand{\cZ}{{\check Z}}

\newcommand{\cCO}{{\mathcal A}(\emptyset)}

\newcommand{\Gr}{\text{\rm Grad}}

\newcommand{\al}{\alpha}

\newcommand{\CS}{\operatorname{CS}}

\theoremstyle{plain}

\newtheorem{theorem}{Theorem}
\newtheorem{proposition}{Proposition}[section]

\newtheorem{corollary}[proposition]{Corollary}

\theoremstyle{definition}

\theoremstyle{remark}

\newtheorem{remark}[proposition]{Remark}

\title[Finite type invariants]{Finite type invariants of 3-manifolds}
\author{Thang T. Q. L\^e}
\begin{document}
\maketitle

\newcommand{\psdraw}[2]
         {\begin{array}{c} \hspace{-1.3mm}
         \raisebox{-4pt}{\psfig{figure=#1.eps,width=#2}}
         \hspace{-1.9mm}\end{array}}

\begin{abstract}
This is a survey article on  finite type invariants of 3-manifolds
written for the \textit{Encyclopedia of Mathematical Physics} to
be published by Elsevier.
\end{abstract}

\section{Introduction}\label{intr}
\subsection{Physics background and motivation}
Suppose $G$ is a semi-simple compact Lie group
$G$  and $M$ a closed oriented 3-manifold. Witten \cite{Witten}
defined quantum invariants by the path integral over all
$G$-connections $A$:

$$Z(M,G;k) := \int \exp(\sqrt{-1}\, k \, \CS(A))\, {\mathcal D} A,$$
where  $k$ is an integer and $\CS(A)$ is the Chern-Simons
functional,

$$ \CS(A) = \frac{1}{4\pi}\int_M Tr(A \wedge dA + \frac{2}{3}
A^3).$$

The path integral is not mathematically rigorous. According to the
stationary phase approximation in quantum field theory, in the
limit $k\to \infty$ the path integral decomposes as a sum of
contributions from the flat connections:

$$ Z(M,G;k) \sim \sum_{\text{flat connections }f} Z^{(f)}(M,K;k) \qquad \text{as $k\to \infty$}.$$

Each contribution is $\exp(2 \pi \sqrt {-1}\, k \,\CS(f))$ times a
power series in $1/k$. The contribution from the trivial
connection is important, especially for rational homology
3-spheres, and the coefficients of the powers $(1/k)^n$,
calculated using $(n+1)$-loop Feynman diagrams by quantum field
theory techniques, are known as perturbative invariants.

\subsection{Mathematical theories} A mathematically rigorous
theory of quantum invariants $Z(M,G;k)$ was pioneered by
Reshetikhin and Turaev in 1990, see \cite{Turaev}. A
number-theoretical expansion of the quantum invariants into power
series that should correspond to the perturbative invariants was
given by Ohtsuki (in the case of $sl_2$, and general simple Lie
algebras by the author) in 1994 that lead him to introducing
finite type invariant theory for 3-manifolds. A universal
perturbative invariant was constructed by Le-Murakami-Ohtsuki (the
LMO invariant) in 1995; it is universal for both finite type
invariants and quantum invariants, at least for homology
3-spheres. Rozansky in 1996 defined perturbative invariants using
Gaussian integral, very close in the spirit to the original
physics point of view. Later Habiro (for $sl_2$ and Habiro and the
author for all simple Lie algebras) found a finer expansion of
quantum invariants, known as the cyclotomic expansion, but no
physics origin is known for the cyclotomic expansion. The
cyclotomic expansion helps to show that the LMO invariant
dominates all quantum invariants for homology 3-spheres.

The purpose of this article is to give an overview of the
mathematical theory of finite type and perturbative invariants of
3-manifolds.

\subsection{Conventions and notations}
All vector spaces are assumed to be over the ground field $\Q$ of
rational numbers, unless otherwise stated. For a graded space $A$,
let $\grad_n A$ be the subspace of grading $n$ and $\grad_{\le
n}A$ the subspace of grading less than or equal to $n$. For $x\in
A$ let $\grad_nx$ and $\grad_{\le n}x$ be the projections of $x$
onto respectively $\grad_n A$ and $\grad_{\le n}A$.

All 3-manifolds are supposed to be closed and oriented. A
3-manifold $M$ is an integral homology 3-sphere (\ZHS) if
$H_1(M,\Z)=0$; it is a rational homology 3-sphere (\QHS) if
$H_1(M,\Q)=0$. For a framed link $L$ in a 3-manifold $M$ denote
$M_L$ the 3-manifold obtained from $M$ by surgery along $L$, see
for example \cite{Turaev}.

%\subsection{Thanks}

\section{Finite type invariants}

After its introduction by Ohtsuki in 1994, the theory of finite
type invariants (FTI) of 3-manifolds  has been developed rapidly
by many authors. Later Goussarov and Habiro independently
introduced clasper calculus, or $Y$-surgery, which provides a
powerful geometric technique  and deep insight in the theory.
$Y$-surgery, corresponding to the commutator in group theory,
naturally gives rise to 3-valent graphs.

\subsection{Generality on finite type invariants}
\subsubsection{Decreasing filtration}
In a theory of FTI, one considers a class of objects, and a
``good" decreasing filtration $\cF_0 \supset \cF_1\supset \cF_2
\supset \dots $ on the vector space $\cF=\cF_0$ spanned by these
objects. An invariant of the objects with values in a vector space
is {\em of order less than or equal to $n$} if its restriction to
$\cF_{n+1}$ is 0; it is of finite type if it is of order less than
or equal to $n$ for some $n$.
 An
invariant has order $n$ if it is of order $\le n$ but not $\le
n-1$. Good here means at least the space of FTI of each order is
finite-dimensional. It is desirable to have an algorithm of
polynomial time to calculate every FTI. Also one wants the set of
FTI invariants to separate the objects (completeness).

The space of invariants of order $\le n$ can be identified with
the dual space of $\cF_0/\cF_{n+1}$; its subspace
$\cF_{n}/\cF_{n+1}$ is isomorphic to the space of invariants of
order $\le n$ modulo the space of invariants of order $\le n-1$.
Informally one can say that $\cF_{n}/\cF_{n+1}$ is more or less
the set of invariants of order $n$.

\subsubsection{Elementary moves, the knot case}
Usually the filtrations are defined using {\em independent
elementary moves}. For the class of knots the elementary move is
given by crossing change. Any two knots can be connected by a
finite sequence of such moves. The idea is if $K,K'\in \cF_{n}$,
the $n$-th term of the filtration, then $K-K'\in \cF_{n+1}$, where
$K'$ is obtained from $K$ by an elementary move. Formal definition
is as follows. Suppose $S$ is a set of double points of a knot
diagram $D$. Let
$$ [D,S] = \sum_{S'\subset S} (-1)^ {\#S'} D_{S'},$$
where the sum is over all subsets $S'$ of $S$, including the empty
set, $D_{S'}$ is the knot obtained by changing the crossing at
every point in $S'$, and $\#S'$ is the number of elements of $S'$.
Then $\cF_n$ is the vector space spanned by all elements of the
form $[D,S]$ with $ \#S=n$. For the knot case, the Kontsevich
integral is an invariant that is universal for all FTI's, see
\cite{BN1}.

\subsection{Ohtsuki's definition of finite type invariants for \ZHS}
An elementary move here is a surgery along a knot: $M \to M_K$,
where $K$ is a framed knot in a \ZHS\ $M$. A collection of moves
corresponds to surgery on a framed link. To always remain in the
class of \ZHS\ we need to restrict ourselves to {\em unit-framed
and algebraically split links}, i.e.  framed links in \ZHS\
 each component of which has framing $\pm1$ and the linking number of
every two components is 0. It is easy to prove that a link $L$ in
a \ZHS\ $M$ is unit-framed and algebraically split if and only if
$M_{L'}$ is a \ZHS\ for every sublink $L'$ of $L$. For a
unit-framed, algebraically split link $L$ in a \ZHS\ $M$ define
$$
[M,L]=\sum_{L'\subset L}(-1)^{|\# L'|}M_{L'},$$
 which is an
element in the vector space  $\cM$ freely spanned by \ZHS.

For  a non-negative integer $n$ let $\cF^{AS}_{n}$ be the subspace
of $\M$ spanned by $[M,L]$ with $\#L=n$. Then the descending
filtration $ \M = \cF^{AS}_{0} \supset \cF^{AS}_{1} \supset
\cF^{AS}_{2} \dots$ defines a theory of FTI on the class of \ZHS.

\begin{theorem}  a) {\rm (Ohtsuki)} The dimension of $\cF_{n}(\cM)$ is
finite for every $n$.

b) {\rm (Garoufalidis-Ohtsuki)} One has $\cF_{3n+1}(\cM)=
\cF_{3n+2}(\cM)=\cF_{3n+3}(\cM)$.
\end{theorem}

 The orders of FTI's in this theory are multiples of 3. The first non-trivial invariant, which is the only
(up to scalar) invariant of degree 3, is the Casson invariant.

\subsection{The Goussarov-Habiro definition}

\subsubsection{$Y$-surgery, or clasper surgery}

\begin{figure}[htb]
$$ \psdraw{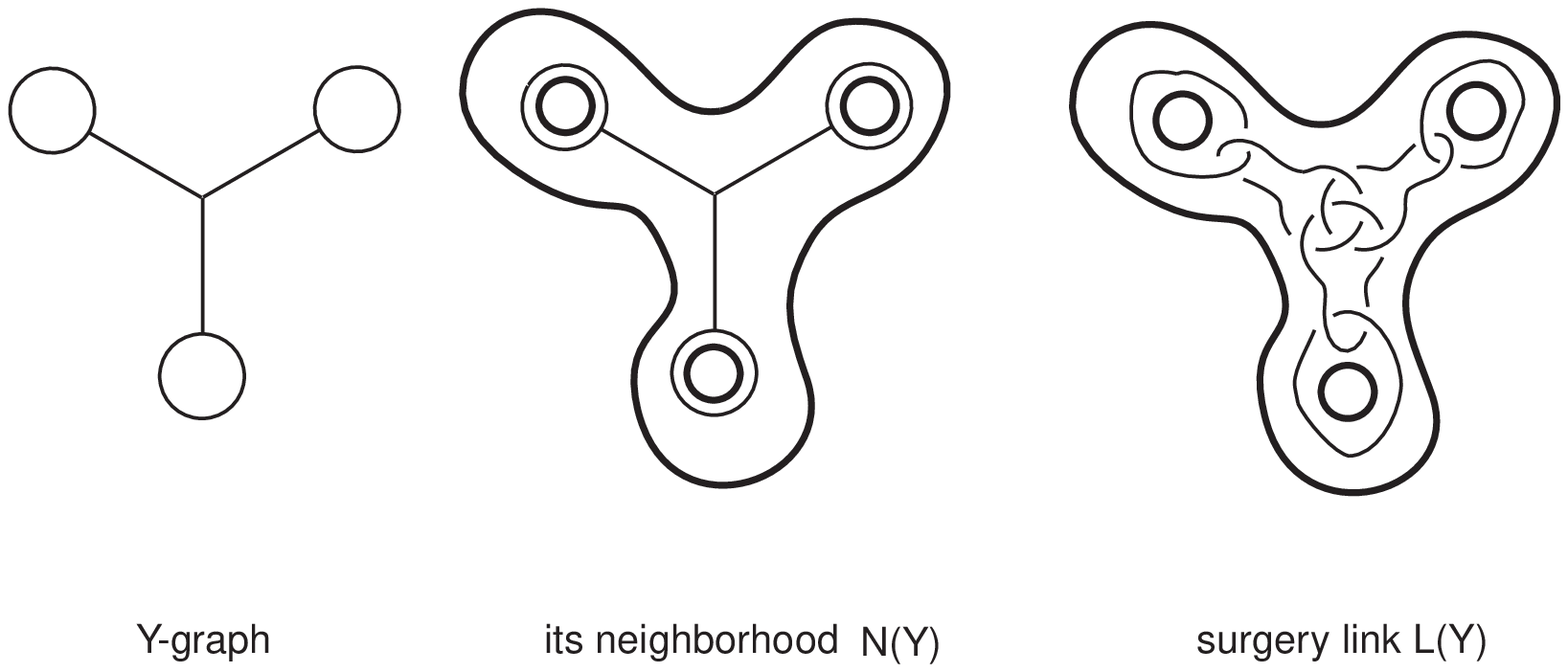}{3.6in} $$ \caption{} \label{13}
\end{figure}

Consider the {\em standard $Y$-graph} $Y$  and a small
neighborhood $N(Y)$ of it in the standard $\R^2$, see Figure
\ref{13}. Denote by $L(Y)$ the 6 component framed link diagram in
$N(Y) \subset \R^3$,  each components of which has framing 0 in
$\R^3$, see Figure \ref{13}.

 A {\em framed $Y$-graph} $C$ in a 3-manifold $M$ is the image of an
embedding of $N(Y)$ into $M$. The surgery of $M$ along the image
of the six-component  link $L(Y)$ is called a {\em $Y$-surgery
along }$C$, denoted by $M_C$. If one of the leaves bounds a disk
in $M$ whose interior is disjoint from the graph, then  $M_C$ is
homeomorphic to $M$.

Matveev in 1987 proved that two 3-manifolds $M$ and $M'$ are
related by a finite sequence of $Y$-surgeries if and only if there
is an isomorphism from $H_1(M,\Z)$ onto $H_1(M',\Z)$ preserving
the linking form on the torsion group. It is natural to partition
the class of 3-manifolds into subclasses of the same $H_1$ and the
same linking form.

\subsubsection{Goussarov-Habiro Filtrations} For a 3-manifold $M$ denote by $\cM(M)$ the vector
space spanned by all 3-manifolds with the same $H_1$ and linking
form.  Define, for a set $S$ of $Y$-graphs in $M$, $[M,S] =
\sum_{S'\subset S} (-1)^{\#S'} M_{S'}$, and $\cF_n^Y\cM(M)$ the
vector space spanned by all $[N,S]$ such that $N$ is in $\cM(M)$
and $\#S=n$. The following theorem of Goussarov and Habiro
(\cite{Gusarov,GGP,Habiro1}) shows that the FTI theory based on
$Y$-surgery is the same as the one of Ohtsuki in the case of \ZHS.

\begin{theorem}
For the case $\cM=\cM(S^3)$, one has $\cF^Y_{2n-1} = \cF^Y_{2n}
=\cF^{AS}_{3n}$.
\end{theorem}

\subsection{The fundamental theorem of finite type invariants of
\ZHS}
\subsubsection{Jacobi diagrams} A {\em closed Jacobi diagram} is a vertex oriented trivalent graph, i.e.
a graph for which the degree of each vertex equal to 3 and a
cyclic order of the 3 half-edges at every vertex is fixed. Here
multiple edges and self-loops are allowed. In pictures, the
orientation at a vertex is the clockwise orientation, unless
otherwise stated. The {\em degree} of Jacobi diagram is half the
number of its vertices.

\begin{figure}[htb]
$$ \psdraw{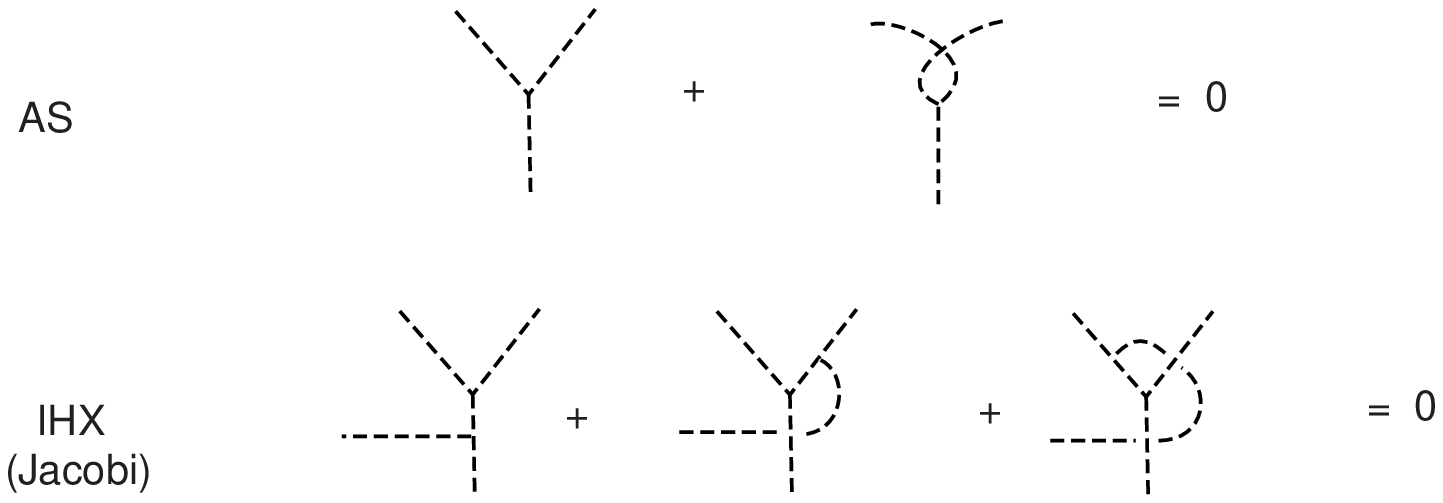}{3.6in} $$ \caption{} \label{AS}
\end{figure}

Let $\grad_n\cA(\emptyset), n\ge 0$, be the vector space spanned
by all closed Jacobi diagrams of degree $n$, modulo the
anti-symmetry (AS) and Jacobi (IHX) relations, see Figure
\ref{AS}.

\subsubsection{The universal weight map $W$}\label{W}
\begin{figure}[htb]
$$ \psdraw{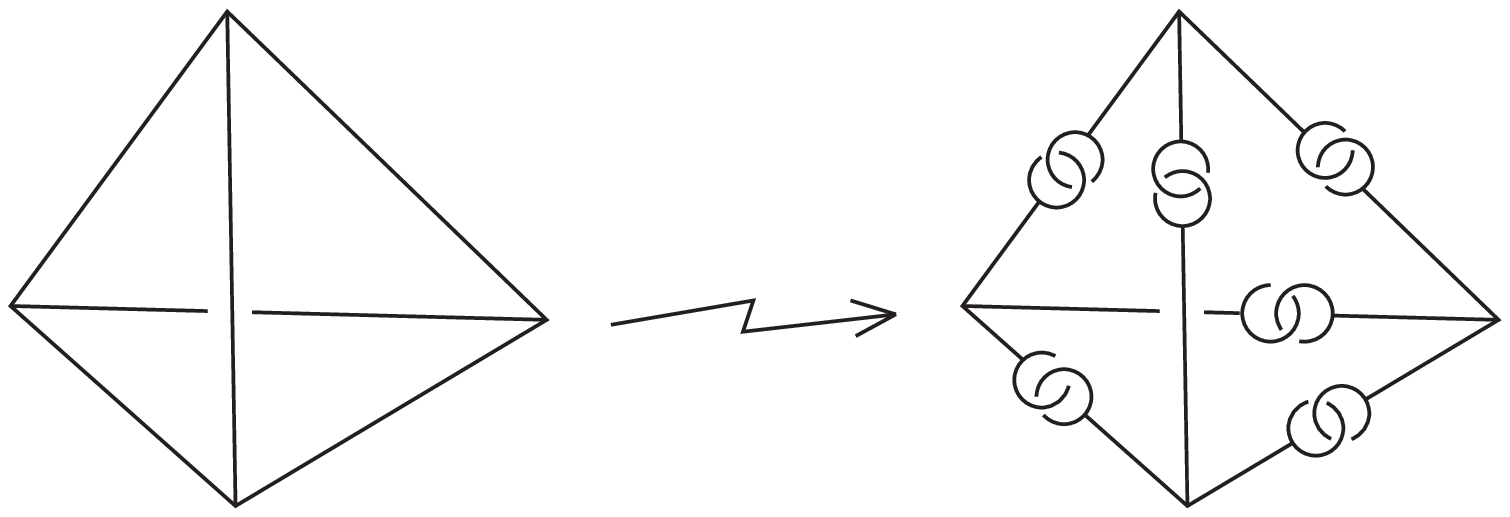}{3.6 in} $$ \caption{} \label{15}
\end{figure}

 Suppose $D$ is a closed Jacobi diagram of degree
$n$.  Embedding $D$ into $\R^3\subset S^3$ arbitrarily and then
projecting down onto $\R^2$ in general position, one can describe
$D$ by a diagram, with over/under-crossing information at every
double point just like in the case of a link diagram. We can
assume that the orientation at every vertex of $D$ is given by
clockwise cyclic order. From the image of $D$ construct a set $G$
of $2n$ $Y$-graphs as in Figure \ref{15}. Here only the cores of a
$Y$-graph are drawn, with the convention that each framed
$Y$-graph is a small neighborhood of its core in $\R^2$.

If $G'$ is a proper subset of $G$, then in $G'$ there is a
$Y$-graph one of the leave of which bounds a disk, hence
$S^3_{G'}=S^3$. Thus $W(D) := [S^3,G]= S^3_{G}-S^3$. By
definition, $W(D) \in \cF^Y_{2n}$; it might depend on the
embedding of $D$ into $\R^3$, but one can show that $W(D)$ is
well-defined in $\cF^Y_{2n}/\cF^Y_{2n+1}$. The map $W$ was first
constructed by Garoufalidis and Ohtsuki in the framework of
$\cF^{AS}$.

\subsubsection{Fundamental theorem}

\begin{theorem} \cite{LMO,Le3} The map $W$ descends to a well-defined linear map $W: \grad_n\cA(\emptyset) \to
\cF^Y_{2n}/\cF^Y_{2n+1}$ and moreover,  is an isomorphism of
vector spaces $\grad_n\cA(\emptyset)$ and
$\cF^Y_{2n}/\cF^Y_{2n+1}$, for $\M=\M(S^3)$.
\end{theorem}

The theorem essentially says that the set of invariants of degree
$2n$ is dual to the space of closed Jacobi diagram
$\grad_n\cA(\emptyset)$. The proof is based on the
Le-Murakami-Ohtsuki invariant (see \S \ref{LMO}).

A $\Q$-valued invariant $I$ of order $\le 2n$ restricts to a
linear map from $\cF_{2n}/\cF_{2n+1}$ to $\Q$. The composition of
$I$ and $W$ is a functional on $\grad_n\cA(\emptyset)$ called the
{\em weight system} of $I$. The theorem shows that every linear
functional on $\grad_{\le n}\cA(\emptyset)$ is the weight of an
invariant of order $\le 2n$.

\subsubsection{Relation to knot invariants} Under the map  that sends an (unframed) knot $K\subset S^3$ to the
\ZHS\
 obtained by surgery along $K$ with framing 1, an
invariant of degree $\le 2n$ (in the $\cF^Y$ theory) of \ZHS\
pulls back to an invariant of order $\le 2n$ of knots. This was
conjectured by Garoufalidis and proved by Habegger.

%\subsubsection{The size of $\grad_n\cA(\emptyset)$}

\subsubsection{Other classes of rational homology
3-spheres} Actually, the theorem was first proved in the framework
of $\cF^{AS}$.  Clasper surgery theory allows Habiro
\cite{Habiro1} to generalize the fundamental theorem to \QHS: For
$M$ a \QHS, the universal weight map $W: \grad_n\cA(\emptyset) \to
\cF_{2n}\cM(M)/\cF_{2n+1}\cM(M)$, defined similarly as in the case
of \ZHS, is an isomorphism, and $\cF_{2n-1}\cM(M)=
\cF_{2n}\cM(M)$.
\subsubsection{Other filtrations and approaches} Other equivalent filtrations were introduced (and compared) by
Garoufalidis, Garoufalidis-Levine and
Garoufalidis-Goussarov-Polyak \cite{GL,GGP}. Of importance is the
one using subgroups of mapping class groups in \cite{GL}. A theory
of $n$-equivalence was constructed by Goussarov and Habiro that
encompasses many geometric aspects of FTI of 3-manifolds
\cite{Habiro1,Gusarov}. Cochran and Melvin \cite{CM} extended the
original Ohtsuki definition to manifolds with homology, using
algebraically split links, but the filtrations are different from
those of Goussarov-Habiro.

\section{The Le-Murakami-Ohtsuki invariant} \label{LMO}

\subsection{Jacobi diagrams}
\begin{figure}[htb]
$$ \psdraw{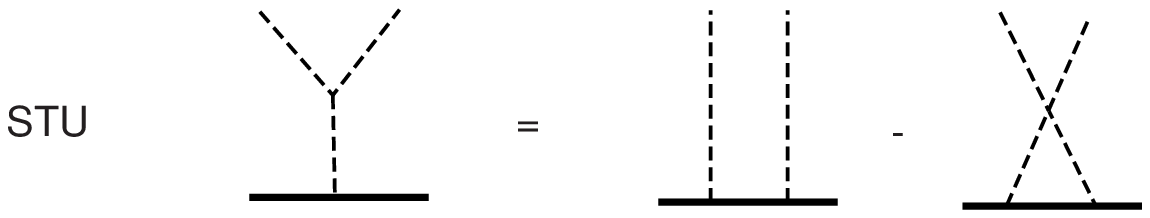}{3in} $$ \caption{} \label{STU}
\end{figure}

 An {\em open Jacobi diagram}  is a
vertex-oriented uni-trivalent graph, i.e., a graph with univalent
and trivalent vertices together with a cyclic ordering of the
edges incident to the trivalent vertices.   A univalent vertex is
also called {\em a leg}. The {\em degree} of an open Jacobi
diagram is half the number of vertices (trivalent and univalent).
A {\em Jacobi diagram based on $X$}, a compact oriented
1-manifold, is a graph $D$ together with a decomposition $D = X
\cup \Gamma$, such that $D$ is the result of gluing all the  legs
of an open Jacobi diagram $\Gamma$ to distinct interior points of
$X$.  The {\em degree} of $D$, by definition, is the degree of
$\Gamma$. In picture  $X$ is depicted by bold lines. Let
$\cA^f(X)$ be the space of Jacobi diagrams based on $X$ modulo the
usual anti-symmetry, Jacobi and the new STU relations (see Figure
\ref{STU}). The completion of $\cA^f(X)$ with respect to degree is
denoted by $\cA(X)$.

When $X$ is a set of $m$ ordered oriented intervals, denote
$\cA(X)$ by $\cP_m$, which has a natural algebra structure where
the product $DD'$ of 2 Jacobi diagrams  is defined by stacking $D$
on top of $D'$ (concatenating the corresponding oriented
intervals). When $X$ is a set of $m$ ordered oriented circles,
denote $\cA(X)$ by $\cA_m$. By identifying the 2 end points of
each interval one gets a map $pr: \cP_m\to \cA_m$,  which is an
isomorphism if $m=1$, see \cite{BN1}.

For $x\in \cA_m$ and $y\in \cA_1$, the connected sum is defined by
$x\#_m y := pr ((pr^{-1}x) \, (pr^{-1}y)^{\otimes m})$, where
$(pr^{-1}y)^{\otimes m}$ is the element in $\cP_m$ with $pr^{-1}y$
on each oriented interval.

\subsubsection{Symmetrization maps}Let $\cB_m$ be the vector space spanned by open Jacobi diagrams
whose legs are labelled by elements of $\{1,2,\dots,m\}$, modulo
the  anti-symmetry and Jacobi relations. One can define an analog
of the Poincare-Birkhoff-Witt
 isomorphism $\chi: \cB_m \to \cP_m$ as follows. For a
diagram $D$, $\chi(D)$ is obtained  by taking the average over all
possible ways of ordering the legs labelled by $j$ and attach them
to the $j$-th oriented interval. It is known that $\chi$ is a
vector space isomorphism \cite{BN1}.

\subsection{The framed Kontsevich integral of links} For an $m$-component framed
link $L\subset \R^3$, the (framed version of the) Kontsevich
integral $Z(L)$ is an invariant  taking values in $\cA_m$, see for
example \cite{Oht}. Let $\nu := Z(K)$, when $K$ is the unknot with
framing 0, and  $\cZ(L):= Z(L) \#_m \nu$. An explicit formula for
$\nu$ is given  in \cite{BLT}.

\subsection{Removing solid loops: the maps $\iota_n$}\label{iota}

Suppose $x\in \cB_m$ is an open Jacobi diagram with legs labelled
by $\{1,\dots,m\}$. If the number of vertices of any label is
different from $2n$, or if the degree of $D$ is greater than
$(m+1)n$, put $\iota_n(D)=0$. Otherwise, partitioning the $2n$
vertices of each label into $n$ pairs and identifying points in
each pair, from $x$ we get a trivalent graph which may contain
some isolated loops (no vertices) and which depends on the
partition. Replacing each isolated loop by a factor $-2n$, and
summing up over all partitions, we get $\iota_n(D)\in \grad_{\le
n}\cA(\emptyset)$.

For $x\in \cA_m$, choose $y\in \cP_m$ such that $pr(y)=x$. Using
the isomorphism $\chi$ we pull back $\chi^{-1}y\in \cB_m$. Define
$\iota_n(x):= \iota_n(\chi^{-1}y)$. One can prove that
$\iota_n(x)$ does not depend on the choice of the preimage $y$ of
$x$. Note that $\iota_n$ lowers the degree by $nm$.

\subsection{Definition of the  Le-Murakami-Ohtsuki invariant $Z^{LMO}$}

In  $\cA(\emptyset):=\prod_{n=0}^\infty \grad_n\cA(\emptyset)$ let
the product of 2 Jacobi diagrams be their disjoint union. Also
define the co-product $\Delta(D) = 1\otimes D+ D\otimes 1$ for $D$
a {\em connected} Jacobi diagram. Then $\cA(\emptyset)$ is a
commutative co-commutative graded Hopf algebra.

For the unknot $U_\pm$  with framing $\pm1$, one has $
\iota_n(\cZ(U_\pm))= (\mp 1)^n+(\text{terms of degree} \ge 1)$,
hence their inverses exist.
 Suppose the linking
matrix of an oriented framed link $L\subset \R^3$ has $\sigma_+$
positive eigenvalues and $\sigma_-$ negative eigenvalues.  Define
\begin{equation}
\Omega_n(L)=\frac{\iota_n(\cZ(L))} {
(\iota_n(\cZ(U_+)))^{\sigma_+} (\iota_n(\cZ(U_-))) ^{\sigma_-}}\in
\Gr_{\le n}(\cCO).\end{equation}

\begin{theorem}{\rm (\cite{LMO})}  $\Omega_n(L)$ is an invariant of the
3-manifold $M= S^3_L$.
\end{theorem}

 We can combine all the $\Omega_n$ to get a better invariant:
$$Z^{LMO}(M):=1+\Gr_1(\Omega_1(M))+\dots+\Gr_n(\Omega_n(M))+\dots
\in \cCO.$$

For $M$ a \QHS, we also define   $$\hat
Z^{LMO}(M):=1+\frac{\Gr_1(\Omega_1(M))}{d(M)}+\dots+
\frac{\Gr_n(\Omega_n(M))} {d(M)^n}+\dots,$$ where
 $d(M)$ is the
 cardinality of
$H_1(M,{\Bbb Z})$.

\begin{proposition}{\rm (\cite{LMO})}\label{grouplike}\quad
Both $Z^{LMO}(M)$ and $\hat Z^{LMO} (M)$ (when defined) are
group-like elements, i.e.
\begin{eqnarray*}\Delta ( Z^{LMO}
(M))&=&Z^{LMO}(M)\otimes
 Z^{LMO}(M),\\
 \Delta (\hat
 Z^{LMO}
(M))&=&\hat Z^{LMO}(M)\otimes
 \hat Z^{LMO}(M).
 \end{eqnarray*}

 Moreover, $\hat Z^{LMO}(M_1\# M_2)= \hat Z^{LMO}(M_1) \times \hat
Z^{LMO}(M_2)$.
 \end{proposition}

\subsection{Universality properties of the LMO invariant}
Let us restrict ourselves to the case of \ZHS.

\begin{theorem} \cite{Le3} The less than or equal to $n$ degree part
$\grad_{\le n} Z^{LMO}$ is an invariant of degree $2n$. Any
invariant of degree $\le 2n$ is a composition  $w(\grad_{\le n}
Z^{LMO})$, where $w: \grad_{\le n}\cA(\emptyset) \to \Q$ is a
linear map.
\end{theorem}

Clasper calculus (or $Y$-surgery) theory allows Habiro to extend
the theorem to rational homology 3-spheres.

\subsection{The Arhus integral}

The Arhus integral (circa 1998) of Bar-Natan, Garoufalidis,
Rozansky and Thurston, based on a theory of formal integration,
calculates the LMO invariant  of rational homology 3-spheres. The
formal integration theory has a conceptual flavor and helps to
relate the LMO invariant to perturbative expansions of quantum
invariants. We give here the definition for the case when one does
surgery on a knot $K$ with {\em non-zero framing} $b$. The link
case is similar, see \cite{BGRT}.

When $K$ is a knot, $\cZ(K)$ is  an element of $\cA_1\equiv
\cP_1\equiv \cB_1$. Note that $\cB_1$ is an algebra where the
product is the disjoint union $\sqcup$. Since the framing is $b$,
one has
$$ \cZ(K) = \exp_{\sqcup}(b\,w_1/2) \sqcup Y,$$ where $w_1$ is the
``dashed interval" (the only connected open Jacobi diagram without
trivalent vertex), and $Y$ is an element in $\cB$ every term of
which must have at least one trivalent vertex. For uni-trivalent
graphs $C,D\in \cB_1$ let $$\langle C,D\rangle =
\begin{cases} 0 \quad \text{if the numbers of legs of $C,D$
are different}\\ \text {sum of all ways to glue legs of $C$ and
$D$ together}
\end{cases}
$$

One defines $\int^{FG} \cZ(K) := \langle \exp_{\sqcup} (- w_1/2b),
Y\rangle$. Then $$\int^{FG} \cZ(K)= \sum_{n=0}^\infty
\frac{\grad_n(\iota_n\,\cZ(K))}{(-b)^n}.$$

And hence

$$\hat Z^{LMO}(S^3_K) = \frac{\int^{FG} \cZ(K)}{\int^{FG}
\cZ(U_{\mathrm{sign } (b)})}.$$

\subsection{Other approaches} Another construction of a universal
perturbative invariant based on integrations over configuration
spaces, closer to the original physics approach but harder to
calculate because of the lack of a surgery formula, was developed
by  Axelrod and Singer, Kontsevich, Bott and Cattaneo, Kuperberg
and Thurston, see \cite{AS, BC}.

\section{Quantum invariants and perturbative expansion}
Fix a simple (complex) Lie algebra $\fg$ of finite dimension.
Using the quantized enveloping algebra of $\fg$ one can define
quantum link  and 3-manifold invariants. We recall here the
definition, adapted for the case of root lattices (projective
group case).

Here our $q$ is equal to $q^2$ in the text book \cite{Jantzen}.
Fix a root system of $\fg$. Let $X,X_+,Y$ denotes respectively the
weight lattice, the set of dominant weights, and the root lattice.
We normalize the invariant scalar product in the real vector space
of the weight lattice so that the length of any short root is
$\sqrt2$.

\subsection{Quantum link invariants}
Suppose $L$ is a framed oriented link with $m$ ordered components,
then the quantum invariant $J_L(\lambda_1,\dots,\lambda_m)$ is a
Laurent polynomial in  $q^{1/2D}$, where
$\lambda_1,\dots,\lambda_m$ are dominant weights, standing for the
simple $\fg$-modules of highest weight
$\lambda_1,\dots,\lambda_m$, and $D$ is the determinant of the
Cartan matrix of $\fg$, see for example \cite{Turaev,Le1}. The
Jones polynomial is the case when $\fg=sl_2$ and all the
$\lambda_i$'s are the highest weight of the fundamental
representation. For the unknot $U$ with 0 framing one has (here
$\rho$ is the half-sum of all positive roots)

\begin{equation*}
J_U(\lambda) = \prod_{\text {positive roots } \al}\frac{
q^{(\lambda+\rho | \al)/2} - q^{-(\lambda +\rho| \al)/2}}{q^{(\rho
| \al)/2} - q^{-(\rho | \al)/2}} \label{unknot1}
\end{equation*}
 We will also
use another normalization of the quantum invariant:

$$Q_L(\lambda_1,\dots,\lambda_m):= J_L(\lambda_1,\dots,\lambda_m) \times \prod_{j=1}^m J_{U}(\lambda_j).$$

This definition is good only for $\lambda_j\in X_+$. Note that
each $\lambda\in X$ is either fixed by an element of the Weyl
group under the dot action (see \cite{Humphreys}) or can be moved
to $X_+$ by the dot action. We define
$Q_L(\lambda_1,\dots,\lambda_m)$ for arbitrary $\lambda_j\in X$ by
requiring that $Q_L(\lambda_1,\dots,\lambda_m)=0$ if one of the
$\lambda_j$'s is fixed by an element of the Weyl group, and that
$Q_L(\lambda_1,\dots,\lambda_m)$ is component-wise invariant under
the dot action of the Weyl group, i.e. for every $w_1,\dots,w_m$
in the Weyl group,

$$ Q_L(w_1\cdot \lambda_1,\dots,w_m\cdot \lambda_m) = Q_L(\lambda_1,\dots,\lambda_m).$$

\begin{proposition} {\rm (\cite{Le1})} Suppose $\lambda_1,\dots,\lambda_m$ are in the
root lattice $Y$.

a) {\rm (Integrality)} Then $Q_L(\lambda_1,\dots\lambda_m) \in
\Z[q^{\pm1}]$, (no fractional power).

b) {\rm (Periodicity)} If $q$ is an $r$-th root of 1, then
$Q_L(\lambda_1,\dots,\lambda_m)$ is invariant under the action of
the lattice group $rY$, i.e. for $y_1,\dots,y_m \in Y$,
$Q_L(\lambda_1,\dots,\lambda_m) =
Q_L(\lambda_1+ry_1,\dots,\lambda_m+ry_m)$.
\end{proposition}

\subsection{Quantum 3-manifold invariants} Although the  infinite sum
$\displaystyle{\sum_{\lambda_j \in Y}
Q_L(\lambda_1,\dots,\lambda_m)}$ does not have a meaning,
heuristic ideas show that it is invariant under the second Kirby
move, and hence almost defines a 3-manifold invariant. The problem
is to regularize the infinite sum. One solution is based on the
fact that at $r$-th roots of unity,
$Q_L(\lambda_1,\dots,\lambda_m)$ is periodic, so we should use the
sum with $\lambda_j$'s running over a fundamental set $P_r$ of the
action of $rY$, where

$$ P_r := \{ x= c_1 \al_1 +\dots + c_\ell \al_\ell  \mid  0\le c_1, \dots,
c_\ell < r\}.$$

Here $\al_1,\dots\al_\ell$ are basis roots. For a root $\xi$ of
 unity of order $r$, let

 $$ F_L(\xi) =  \sum_{\lambda_j\in (P_r\cap Y)}
Q_L({\lambda_1},\dots,{\lambda_m})|_{q=\xi}.$$

For an oriented framed link $L$ let $\sigma_+$ and $\sigma_-$ be
respectively the number of positive and negative eigenvalues of
the linking matrix of $L$. Let $U_\pm$ be the unknot with framing
$\pm 1$. If $F_{U_\pm}(\xi) \neq 0$, define

\begin{equation*}
\tau_L(\xi) := \frac{F_L(\xi)}{(F_{U_{+}}(\xi))^{\sigma_+}
(F_{U_-}(\xi))^{\sigma_-}} \label{definition3}.\end{equation*}

Recall that $D$ is the determinant of the Cartan matrix. Let $d$
be the maximum of the absolute values of entries of the Cartan
matrix outside the diagonal.
\begin{theorem} \cite{Le2} a) If the order $r$ of $\xi$ is co-prime with $dD$,
then $F_{U_\pm}(\xi) \neq 0$.

b) If $F_{U_\pm}(\xi) \neq 0$ then $\tau^{P\fg}_M(\xi):=
\tau_L(\xi)$ is an invariant of the 3-manifold $M=S^{3}_L$.
\end{theorem}

\begin{remark} The version presented here corresponds to projective groups. It was defined by Kirby and Melvin for
 $sl_2$, Kohno and Takata for $sl_n$, and for arbitrary simple Lie algebra
by the author \cite{Le2}. When $r$ is co-prime with $d D$, there
is also an associated modular category that generates a
Topological Quantum Field Theory. In most texts in literature, say
\cite{Kirr,Turaev}, another version $\tau^\fg$ was defined. The
reason we choose $\tau^{P\fg}$ is it has nice integrality and
eventually a perturbative expansion. For relations between the
version $\tau^{P\fg}$ and the usual $\tau^\fg$, see \cite{Le2}.
\end{remark}

\subsubsection{Examples}
 When $M$ is the Poincare sphere and $\fg= sl_2$,

$$ \tau^{Psl_2}_M(q) = \frac{1}{1-q}\sum_{n=0}^\infty q^n (1-q^{n+1})(1-q^{n+2}) \dots
(1-q^{2n+1}).
$$

Here $q$ is  a root of unity, and the sum is easily seen to be
finite.

\subsubsection{Integrality} The following theorem was proved for  $\fg =sl_2$ by H.
Murakami \cite{Murakami}
 and $\fg =sl_n$ by Takata-Yokota
and Masbaum-Wenzl (using ideas of J. Roberts) and for arbitrary
simple Lie algebras by the author in \cite{Le2}.

\begin{theorem} Suppose the order $r$ of $\xi$ is a prime big enough,
then $\tau_M^{P\fg}(\xi)$ is in $\Z[\xi]=\Z[\exp(2\pi i/r)]$.
\label{integ3}
\end{theorem}

\subsection{Perturbative expansion}
 Unlike the link case, quantum 3-manifold
invariants can be defined only at certain roots of unity. In
general,  there is no analytic extension of the function
$\tau^{P\fg}_M$ around $q=1$. In perturbative theory, we want to
expand the function $\tau^\fg_M$ around $q=1$ into power series.
For \QHS, Ohtsuki (for $\fg=sl_2$) and then the author (for all
other simple Lie algebras) showed that there is a
number-theoretical expansion of $\tau^{P\fg}_M$ around $q=1$ in
the following sense.

 Suppose $r$ is a big enough prime, and
$\xi = \exp (2 \pi i/r)$. By the integrality (Theorem
\ref{integ3}),

$$\tau^{P\fg}_M(\xi) \in \Z[\xi] =  \Z[q] /(1 +q + q^2 + \dots + q^{r-1}).$$

Choose a representative $f(q)\in \Z[q]$ of $\tau^{P\fg}_M(\xi)$.
Formally substitute $q=(q-1)+1$ in $f(q)$:

$$f(q) = c_{r,0} + c_{r,1} (q-1) + \dots+  c_{r,n-2} (q-1)^{n-2} $$

The integers $c_{r,n}$ depend on $r$ and the representative
$f(q)$.  It is easy to see that   $c_{r,n} \pmod r$ does not
depend on the representative $f(q)$ and hence is an  invariant of
\QHS. The dependence on $r$ is a big drawback. The theorem below
says that there is a {\em rational} number $c_n$, not depending on
$r$, such that $c_{r,n} \pmod r$ is the reduction of either $c_n$
or $-c_n$ modulo $r$, for sufficiently large prime $r$. It is easy
to see that if such $c_n$ exists, it must be unique.  Let $s$ be
the number of positive roots of $\fg$. Recall that $\ell$ is the
rank of $\fg$.

\begin{theorem} For every\  {\QHS}\ $M$, there is a sequence of numbers $c_n \in \Z[\frac{1}{(2n+ 2s)!
|H_1(M,\Z)|}]$, such that for sufficiently large prime $r$

$$ c_{r,n}  \equiv \left( \frac{|H_1(M,\Z)|}{r}\right)^\ell \, c_n \pmod r,$$ where
$\left(\frac {|H_1(M,\Z)|}{r}\right)=\pm 1 $ is the Legendre
symbol. Moreover, $c_n$ is an invariant of order $\le 2n$.
\end{theorem}

The series  ${\frak t}^{P\fg}_M(q-1) := \sum_{n=0}^\infty c_n
(q-1)^n$, called the Ohtsuki series, can be considered as the
perturbative expansion of the function $\tau_M^{P\fg}$ at $q=1$.
For actual calculation of ${\frak t}^{P\fg}_M(q-1)$ see
\cite{Le2,Oht,Roz}.

\subsubsection{Recovery from the LMO invariant}
It is known that for any metrized Lie algebra $\frakg$, there is a
linear  map $W_{\frakg}:\grad_n\cA(\emptyset) \to \Q$, see
\cite{BN1}.

\begin{theorem} One has
$$\sum_{n=0}^\infty W_\fg(\grad_nZ^{LMO})\, h^n = {\frak t}^{P\fg}_M(q-1)|_{q=e^h}.$$
\end{theorem}

This shows that the Ohtsuki series ${\frak t}^{P\fg}_M(q-1)$ can
be recovered from, and hence totally determined by, the LMO
invariant. The theorem was proved by Ohtsuki for $sl_2$. For other
simple Lie algebras  the theorem follows from the Arhus integral,
see \cite{BGRT,Oht}.

\subsection{Rozansky's Gaussian integral}
Rozansky gave a definition of the Ohtsuki series using formal
Gaussian integral in the important work \cite{Roz}. The work is
only for $sl_2$, but can be generalized to other Lie algebras; it
is closer to the original physics ideas of perturbative
invariants.

\section{Cyclotomic Expansion}
\subsection{The Habiro ring}

Let us define the
\newcommand{\Ha}{\widehat{\Z[q]}}
Habiro ring $\Ha$ by

$$\Ha := \lim_{\leftarrow n}\Z[q]/((1-q)(1-q^2)\dots(1-q^n)).$$

Habiro \cite{Habiro2} called it the cyclotomic completion of
$\Z[q]$. Formally, $\Ha$ is the set of all series of the form

$$f(q) = \sum_{n=0}^\infty f_n(q) \, (1-q)(1-q^2)\dots(1-q^n), \qquad \text{where } \quad f_n(q) \in \Z[q].$$

 Suppose $U$ is the set of roots of 1. If $\xi\in U$ then
 $(1-\xi)(1-\xi^2) \dots (1-\xi^n)=0$ if $n$ is big enough, hence
one can define $f(\xi)$ for $f\in \Ha$. One can consider every
$f\in \Ha$ as a function with domain $U$. Note that $f(\xi)\in
\Z[\xi]$ is always an algebraic integer.  It turns out $\Ha$ has
remarkable properties, and plays an important role in quantum
topology.

Note that the formal derivative of $(1-q)(1-q^2)\dots(1-q^n)$ is
divisible by $(1-q)(1-q^2)\dots(1-q^k)$ with $k$ the integer part
of $(n-1)/2$. This means every element $f\in \Ha$ has a derivative
$f'\in \Ha$, and hence derivatives of all orders in $\Ha$. One can
then associate to $f\in \Ha$ its Taylor series at a root $\xi$ of
1:

$$T_\xi(f) := \sum_{n=0}^\infty
\frac{f^{(n)}(\xi)}{n!}(q-\xi)^n,$$ which can also be obtained by
noticing that $(1-q)(1-q^2)\dots(1-q^n)$ is divisible by
$(q-\xi)^k$ if $n$ is bigger than $k$ times the order of $\xi$.
Thus one has a map  $T_\xi :\Ha \to \Z[\xi][[q-\xi]] $.

\begin{theorem}\cite{Habiro2} a) For each root of unity $\xi$, the map $T_\xi$ is
injective, i.e. a function in $\Ha$ is determined by its Taylor
expansion at a point in the domain $U$.

b) If $f(\xi)=g(\xi)$ at infinitely many roots $\xi$ of prime
power orders, then $f=g$ in $\Ha$.

\end{theorem}

One important consequence is that $\Ha$ is an integral domain,
since we have the embedding $T_1 :\Ha \hookrightarrow \Z[[q-1]]$.

In general the Taylor series $T_1f$ has 0 convergence radius.
However, one can speak about $p$-adic convergence to $f(\xi)$ in
the following sense. Suppose the order $r$ of $\xi$ is a power of
prime, $r=p^k$. Then it's known that $(\xi-1)^n$ is divisible by
$p^m$ if $n > mk$. Hence $T_1f(\xi)$ converges in the $p$-adic
topology, and it can be easily shown that the limit is exactly
$f(\xi)$.

 The above properties suggest to consider $\Ha$ as a class of ``analytic functions" with domain
$U$.

\subsection{Quantum invariants as an element of $\Ha$}

It was proved, by Habiro for $sl_2$ and by Habiro with the author
for general simple Lie algebras, that quantum invariants of \ZHS's
belong to $\Ha$ and thus have remarkable integrality properties:
\begin{theorem} a) For every {\ZHS}\ $M$, there is an invariant
$I^\fg_M\in \Ha$ such that if $\xi$ is a root of unity for which
the quantum invariant $\tau^{P\fg}_M(\xi)$ can be defined, then
$I^\fg_M(\xi) = \tau^{P\fg}_M(\xi)$.

b) The Ohtsuki series is equal to the Taylor series of $I^\fg_M$
at 1.
\end{theorem}

\begin{corollary} Suppose $M$ is a {\ZHS}.

a) For every root of unity $\xi$, the quantum invariant at $\xi$
is an algebraic integer, $\tau^\fg_M(\xi)\in \Z[\xi]$. (No
restriction on the order of $\xi$ is required).

b) The Ohtsuki series ${\frak t}^{P\fg}_M(q-1)$ has integer
coefficients. If $\xi$ is a root of
 order $r=p^k$, where $p$ is prime, then the Ohtsuki series at
 $\xi$
 converges $p$-adically to the quantum invariant at $\xi$.

 c) The quantum invariant $\tau^{P\fg}_M$ is determined by values at infinitely many roots
 of prime power orders and also determined by its
 Ohtsuki series.

 d) The LMO invariant totally determines the quantum invariants $\tau^{P\fg}_M$.
\end{corollary}

 Part (b) was conjectured by R. Lawrence for $sl_2$ and
first proved by Rozansky (also for $sl_2$ \cite{Rozansky}). Part
(d) follows from the fact that the LMO invariant determines the
Ohtsuki series; it exhibits  another universality property of the
LMO invariant.

\end{document}